\documentclass{amsart}

\usepackage{amssymb}
\newtheorem{thm}{Theorem}
\newtheorem{lem}{Lemma}

\begin{document}

\bibliographystyle{plain}

\title[Romi  Shamoyan ]{On an extremal problem in analytic  spaces in tube domains over symmetric cones}

\author[]{Romi  Shamoyan}

\address{Department of Mathematics, Bryansk State Technical University, Bryansk ,241050, Russia}
\email{\rm rshamoyan@gmail.com}

\date{}

\begin{abstract}
New sharp estimates concerning distance function in certain Bergman -type spaces of analytic functions on tube domains over symmetric cones are obtained.
These are first results of this type for tube domains over symmetric cones.
\end{abstract}

\maketitle

\footnotetext[2]{\, Mathematics Subject Classification 2010 Primary 42B15, Secondary 42B30.  Key words
and Phrases: Distance estimates,tube domains}

\section{Introduction and preliminaries}

In this note we obtain sharp distance estimates in spaces of analytic spaces in tube domains over symmetric cones.
 
This line of investigation can be considered as a continuation of previous papers \cite{AS1}, \cite{SM1} and \cite{SM2}. 

These new results are 
contained in the second section of this note. We remark that  for the first time in literature we consider this extremal problem related with distance estimates in spaces of analytic functions on tube domains over symmetric cones. 
The first section  contains  required preliminaries on analysis on symmetric cones.

In one dimensional tubular domain which is  upperhalfspace  $C_{+}$ (see \cite{BBR}) our theorems are not new and they were obtained recently in \cite{SA}.

Moreover arguments in proofs we provided below are similar to those we have in one dimension
and the base of proof is again the so-called Bergman reproducing formula, but in tubular domain over symmetric cone.(see ,for example, \cite{BBR} for this integral representation)
 
We shortly remind the history of this problem.

After the appearance of \cite {Zh} various papers appeared where  arguments which can be seen in \cite{Zh} were extended and modifyied in various directions
\cite {SM1},\cite{SM2},\cite{AS1}. 

In particular in mentioned papers various new results on distances for analytic function spaces in higher dimension (unit ball and polydisk ) were obtained.
Namely new results  for large scales of analytic mixed norm spaces in higher dimension were proved.

Later several new sharp results for harmonic functions of several variables in the unit ball and upperhalfplane of Euclidean space were also obtained 
see ,for example, \cite{AS1} and references there.

We mention separatly \cite{SA} and \cite {SR} where the case of higher dimension was considered in  special cases
of analytic spaces on subframe and new analogues results in the context of bounded strictly pseudoconvex domains with smooth boundary  were also provided.

The classical Bergman representation formula in various forms and in various domains serves as a base in all these papers in proofs of main results.

We would like to note also recently Wen Xu \cite{WX} repeating arguments of Ruhan Zhao in the unit ball obtained results on distances
from Bloch functions to some Mobius invariant function spaces in one and higher dimension in a relativly direct way.

Probably for the first time in literature these extremal problems connected with distances in analytic spaces appeared before in \cite{ACP} and in \cite{An}

where this problem was formulated and certain cases connected with spaces of bounded analytic functions in the unit disk were considered.

These results were mentioned \cite {X1} ,some other results on distance problems in BMOA spaces can be found also in \cite{GZ}.

VArious other  extremal problems in analytic function spaces also were considered  before in various papers, see for example \cite{AL},\cite{KH},\cite{RU},\cite{KS}. 

In those papers other results around this topic and  some applications of certain extremal problems can be found also.

The goal of this note to develop further  some ideas from our recent mentioned papers and present a new sharp theorem in tube domain over symmetric cones.

We note in case of upper halfplane of complex plane $C$ which is a tube domain in one dimension
such results already were obtained previously by author \cite{SA}.
For formulation of our results we will need various standard definitions from the theory of tube domains over symmetric cones.( see \cite{DD},\cite{FaKo},\cite{SE},\cite{BBR})

Let $T_\Omega = V + i\Omega$ be the tube domain over an irreducible symmetric cone $\Omega$ in the
complexification $V^{\mathbb C}$ of an $n$-dimensional Euclidean space $V$. $\mathcal H(T_\Omega)$ denotes the space of
all holomorphic functions on $T_\Omega$. Following the notation of
\cite{FaKo}	and \cite {BBR} we denote the rank of the cone $\Omega$ by $r$ and by $\Delta$ the determinant function on $V$.

Letting $V = \mathbb R^n$, we have as an example of a symmetric cone on $\mathbb R^n$ the Lorentz cone
$\Lambda_n$ which is a rank $2$ cone defined for $n \geq 3$ by
$$\Lambda_n = \{ y \in \mathbb R^n : y_1^2 - \cdots - y_n^2 > 0, y_1 > 0 \}.$$

The determinant function in this case is given by the Lorentz form
$$\Delta(y) = y_1^2 - \cdots - y_n^2.$$ (see for example \cite{BBR})

Let us introduce some convinient notations regarding multi-indicses.

If $t = (t_1, \ldots, t_r)$, then $t^\star = (t_r, \ldots, t_1)$ and, for $a \in \mathbb R$, $t+a =
(t_1 + a, \ldots, t_r + a)$. Also, if $t, k \in \mathbb R^r$, then $t < k$ means $t_j < k_j$ for all
$1 \leq j \leq r$.

We are going to use the following multi-index
$$g_0 = \left( (j-1)\frac{d}{2} \right)_{1 \leq j \leq r}, \;\;\; \mbox{\rm where} \;\;\;\
(r-1) \frac{d}{2} = \frac{n}{r} - 1.$$

For $\tau \in \mathbb R_{+}$ and the associated determinant function $\Delta(x)$ \cite{BBR} we set
\begin{equation}
A^\infty_\tau (T_\Omega) = \left\{ F \in {\mathcal H}(T_\Omega) : \| F \|_{A^\infty_\tau} = \sup_{x+iy \in T_\Omega}
|F(x+iy)| \Delta^\tau (y) < \infty \right\},
\end{equation}

It can be checked that this is a Banach space .
Below we denote by   $\Delta_s$  the generalized power function \cite{FaKo},\cite{BBR}. 

For $1 \leq p, q < + \infty$ and $\nu \in \mathbb R$,and $\nu > \frac{n}{r}-1$ we denote by $A_\nu^{p,q}(T_\Omega)$ the mixed-norm weighted
Bergman space consisting of analytic functions $f$ in $T_\Omega$ such that
$$ \| F \|_{A_\nu^{p,q}} = \left( \int_\Omega \left( \int_V |F(x+iy)|^p d\,x \right)^{q/p} \Delta^\nu (y)
\frac{d\, y}{\Delta (y)^{n/r}} \right)^{1/q} < \infty,$$

This is a Banach space.
 
Replacing above A by L we will get as usual the corresponding larger space of all  measurable functions in tube over symmetric cone with the same quazinorm
(see \cite{FaKo},\cite{SE})

It is known the $A_\nu^{p,q}(T_\Omega)$ space is nontrivial if and only if $\nu > \frac{n}{r}-1 $,(see \cite{DD},\cite{BBR}). 

When $p = q$ we write (see \cite{BBR})
$$A_\nu^{p,q}(T_\Omega) = A^p_\nu(T_\Omega)$$. This is the classical weighted Bergman space with usual modification when $p=\infty$ 

The (weighted) Bergman projection $P_\nu$ is the orthogonal projection from the Hilbert space $L^2_\nu(T_\Omega)$
onto its closed subspace $A^2_\nu(T_\Omega)$ and it is given by the following integral formula (see \cite{BBR})

\begin{equation}\label{bpro}
P_\nu f(z) = C_\nu \int_{T_\Omega} B_\nu (z, w) f(w) d V_\nu (w),
\end{equation}

where $$B_\nu (z, w) = C_\nu \Delta^{-(\nu + \frac{n}{r})}((z-\overline w)/i)$$ is the Bergman reproducing kernel

for $$A^2_\nu(T_\Omega)$$,(see \cite{FaKo},\cite{BBR})

Here we used the notation $$d V_\nu (w) = \Delta^{\nu - \frac{n}{r}}(v) du dv$$.

Below and here we use constantly the following notations $ w = u +iv \in T_\Omega $ and also $ z=x+iy \in T_\Omega $.

Hence for any analytic function from $A^2_\nu(T_\Omega)$ the following integral formula  is valid.(see also \cite{BBR})

\begin{equation}\label{bpro1}
f(z)=C_\nu \int_{T_\Omega} B_\nu (z, w) f(w) d V_\nu (w),
\end{equation}

In this case  sometimes below we say simply that the $f$ function allows Bergman representation via Bergman kernel with $\nu$ index.

Note these assertions have direct copies in simpler cases of analytic function spaces in unit disk,polydisk,unit ball,upperhalfspace $C_{+}$
and in spaces of harmonic functions in the unit ball or upperhalfspace of Euclidean space $R^{n}$
These classical facts are well- known and can be found ,for example ,in \cite{DS} and in some items from references there.

Above and throughout the paper  we write $C$( sometimes with indexes ) to denote   positive  constants which might be different  each time we see them 
(and even in a chain of inequalities), but is independent of the functions or variables being discussed.

In this paper we will also need a pointwise estimate for the Bergman projection of functions in $$L^{p,q}(T_\Omega)$$, defined by
integral formula (\ref{bpro}), when this projection makes sense.Note such estimates in simpler cases of unit disk, unit ball and polydisk are well- known (see \cite{DS})

 Let us first recall the following known basic integrability properties for the determinant function, which appeared already above in definitions.

\begin{lem}
Let $\alpha \in \mathbb C^r$ and $y \in \Omega$.

1) The integral
$$J_\alpha (y) =  \int_{\mathbb R^n} \left| \Delta_{-\alpha} \left( \frac{x+iy}{i} \right) \right| dx$$
converges if and only if $\Re \alpha > g_0^\ast + \frac{n}{r}$. In that case $$J_\alpha (y) = C_\alpha
|\Delta_{-\alpha + n/r} (y)|$$.

2) For any multi-indices $s$ and $\beta$ from $C^{r}$ and $t \in \Omega$ the function $$y \mapsto \Delta_\beta(y+t) \Delta_s(y)$$
belongs to $L^1(\Omega, \frac{dy}{\Delta^{n/r}(y)})$ if and only if $\Re s > g_0$ and $\Re (s+\beta) < - g_0^\ast$. In
that case we have
$$\int_\Omega \Delta_\beta (y+t) \Delta_s(y) \frac{dy}{\Delta^{n/r}(y)} = C_{\beta, s} \Delta_{s+\beta}(t).$$
\end{lem}

We refer to Corollary 2.18 and Corollary 2.19 of \cite{DD} for the proof of the above lemma or \cite{BBR} 

As a corrolary of one dimensional versions of these estimates (see, for example, \cite{SE}Theorem 3.9) we obtain the following vital estimate (A) which we will use in proof of our main result.

$$\int_{T_\Omega }\Delta^{\beta}(y)|B_{\alpha+\beta+\frac{n}{r}}(z,w)|dV(z)\leq C \Delta ^{-\alpha}(v)$$,$\beta>-1$,
$\alpha>\frac{n}{r}-1$,$z=x+iy$,$w=u+iv$ (see \cite{SE})

Let $\tau$ be the set
of all triples $(p, q, \nu)$ such that $1 \leq p, q < \infty$, $\nu > \frac{n}{r}-1$ . 

The following vital pointwise estimate can be found ,for example, in \cite{BBR}.

\begin{lem}\label{Ltau}

Suppose $(p, q, \nu) \in \tau$. Then

\begin{equation}\label{poin}
|P_\nu f(z)| \leq  c_{p,q,r,\nu,n} \Delta^{-\frac{\nu}{q} - \frac{n}{rp}}(\Im z) \| f \|_{A^{p,q}_\nu}.
\end{equation}

\end{lem}

{\it Proof.} This is a consequence of the  lemma formulated above and H\"older's inequality (see \cite{BBR}) $\Box$

\section {New estimates for distances in analytic function spaces in tube domains over symmetric cones}

In this paper we restrict ourselves to  $\Omega$  irreducible symmetric cone in the Euclidean vector space $R^{n}$ of dimension n,endowed with an inner product for which
the cone $\Omega$ is self dual.We denote by $T_\Omega=R^{n}+i\Omega$ the corresponding tube domain in $C^{n}$.

This section is devoted to formulations and proofs of all main results of this paper.As previously in case of analytic functions in 
unit disk,polydisk,unit ball, and upperhalfspace $C_{+}$ and in case of spaces of harmonic functions in Euclidean space \cite {Zh},\cite{SR},\cite{SA},\cite{AS1},\cite{SM1},
\cite{SM2} the role of the Bergman representation formula is crucial in these issues and our proof is heavily based on it.

As it is known a variant of  Bergman representation formula is available also in Bergman- type analytic function 
spaces in tubular domains over symmetric connes and this known fact (see \cite{DD},\cite{FaKo},\cite{BBR},\cite{SE}), which is crucial also in various problems
in analytic function spaces in tubular domains (see \cite{BBR} and various references there) is used also in our proofs below.

The following result can be found in \cite{SE}(section4)

For all $1< p < \infty $ and  $1< q < \infty$ and for all $\frac{n}{r} \leq p_{1}$ ,where $\frac{1}{p_{1}}+\frac{1}{p}=1$ ,and $$ \frac{n}{r}-1 < \nu $$

and for all functions $f$ from $A^{p,q}_\nu$ and for all $ \frac{n}{r}-1 < \alpha $ the Bergman  representation formula with $\alpha$ index 

or with the Bergman kernel $$B_\alpha(z,w)$$ is valid .

We remark this result is a particular case of a more general assertion for analytic mixed norm $A^{p,q}_\nu $ classes (see \cite{SE} )
which means that our main result below admits  also some extentions,even to mixed norm spaces which we defined above,this will be discussed  below at the end and in our next paper which is in preparation.   

We will also need for our proofs the following important fact on integral representations.(see \cite{BBGRS})
Let $\nu>\frac{n}{r}-1$,$\alpha>\frac{n}{r}-1$,then for all functions from $A^\infty_\alpha$ the integral represenations of Bergman with Bergman kernel $$B_{\alpha+\nu}(z,w)$$
(with $\alpha+\nu$ index) is valid.

We note also that by Lemma \ref{Ltau} we have
\begin{equation}
|f(x + iy)| \Delta^{\frac{n}{rp} + \frac{\nu}{q}}(y) \leq c_{p,q,r,\nu} \| f \|_{A^{p,q}_\nu}, \qquad x + iy \in T_\Omega, \quad
(p,q,\nu) \in \tau.
\end{equation}

This means that we have a continuous embedding $A^p_\nu \hookrightarrow A^\infty_{\frac{n}{rp} + \frac{\nu}{p}} $ for
$(p,p,\nu) \in \tau$ and this naturally leads to a problem of estimating $${\rm dist}_{A^\infty_{\frac{n}{rp} + \frac{\nu}{p}}}
(f, A^p_\nu)$$ for a given $f \in A^\infty_{\frac{n}{rp} + \frac{\nu}{p}}$. 

This problem is solved in our next theorem below, 
which is the main result of this section. Let us set, for $f \in {\mathcal H}(T_\Omega)$, $s \in \mathbb R$ and $\epsilon > 0$:

\begin{equation}\label{Efsep}
V_{\epsilon, s}(f) = \left\{ x + iy \in T_\Omega : |f(x+iy)|\Delta^s(y) \geq \epsilon \right\}
\end{equation}
Let  also $w=u+iv \in T_\Omega $,$z=x+iy \in T_\Omega $.We denote by $N_{1}$ and by $N_{2}$ two sets- the first one is $V_{\epsilon,s}(f)$, the other one
is the set of all those points ,which are in tubular domain $T_\Omega$, but not in $N_{1}$

\begin{thm}\label{Td1}
Let $1 < p < \infty$,  $\nu > p(\frac{n}{r}-1)$,$\beta > t+\frac{n}{r}-1$ ,
$t = \frac{1}{p} (\nu + \frac{n}{r})$.Set, for $f \in A^\infty_{\frac{n}{rp} + \frac{\nu}{p}}$:

\begin{equation}\label{ETd1a}
l_1(f) = {\rm dist}_{A^\infty_{\frac{n}{rp} + \frac{\nu}{p}}} (f, A^p_\nu),
\end{equation}
\begin{equation}\label{ETd1b}
l_2(f) = \inf \left\{ \epsilon > 0 : \int_{T_\Omega} \left( \int_{V_{\epsilon, t}(f)} \frac{\Delta^{\beta - t - \frac{n}{r}}(v)
dudv}{\Delta^{\beta + \frac{n}{r}}((z-\overline w) /i )} \right)^p \Delta^{\nu-\frac{n}{r}}(y) dxdy < \infty \right\}.
\end{equation}

Then there is a positive number $ \beta_{0}$   ,so that for all $ \beta > \beta_{0}$   we have $l_1(f) \asymp l_2(f)$.
\end{thm}

{\it Proof.}We will use for our proofs the following observation.\cite{BBGRS}
Let $\nu>\frac{n}{r}-1$,$\tau>\frac{n}{r}-1$,then for all functions from $A^\infty_\tau$ the integral represenations of Bergman with Bergman kernel $$B_(\tau+\nu)(z,w)$$
is valid.  

 We denote below the double integral which appeared in formulation by $G(f)$ and we will show first that  $l_{1}(f)\leq C l_{2}(f) $ .We assume now that $l_{2}(f)$ is finite. 

 We use the Bergman representation formula which we provided above ,namely(\ref{bpro1}) ,and using conditions on parameters we have  the following  equalities .

First we have obviously by remark from which we started this proof that for large enough $\beta$
$$f(z)=C_\nu\int_{T_\Omega}B_\beta(z,w)f(w)dV_\beta(w)= f_{1}(z)+f_{2}(z)$$
$$f_{1}(z)=C_\nu\int_{N_{2}}B_\beta(z,w)f(w)dV_\beta(w)$$,$$f_{2}(z)=C_\nu\int_{N_{1}}B_\beta(z,w)f(w)dV_\beta(w)$$

Then we estimate both functions separatly using lemmas provided above and following some arguments we provided in one dimensional case that is the case of upperhalfspace $C_{+}$ 
\cite{SA}.Here our arguments are sketchy since they are parallel to arguments from \cite{SA}. 
Using definitions of $N_{1}$ and $N_{2}$ above  after some calculations  following arguments from \cite {SA} using estimate (A) we will have immediatly.

$$f_{1}\in A^\infty_{\frac{n}{rp}+\frac{\nu}{p}}$$   and  $$f_{2}\in A^p_\nu$$.We easily note the last inclusion follows directly from the fact that $l_{2}$ is finite.

Moreover it can be easily seen that the norm of $f_{1}$ can be estimated from above by $C\epsilon$,where 
for some positive constant $C$ (\cite{SA}),since obviously $$\sup_{N_{2}}|f(w)|\Delta^{t}(v)\leq \epsilon$$
Note this last fact follows directly from definition of $N_{2}$ set and estimates in lemma above
which leads to the following  inequality which was denoted by us as (A).(see also \cite{SE})

$$\int_{T_\Omega}\Delta ^{-t}(y)|B_{\beta}(z,w)|dV_{\beta}(z)\leq C \Delta ^{-t}(v)$$,$z=x+iy$,$w=u+iv$,

for all $\beta$ so that $\beta>\beta_{0}$ ,for some large enough fixed $\beta_{0}$ which depends on $n$,$r$,$\nu$,$p$ and  
for $t=(\frac{1}{p})(\nu+\frac{n}{r})$ and $\nu>p(\frac{n}{r}-1)$ (see \cite{SE} Theorem 3.9)

This gives immediatly one part of our theorem .Indeed, we have obviously

$$l_{1} \leq C_{2} \|f-f_{2}\|_{A^\infty_t} =C_{3} \|f_{1}\|_{A^\infty_t} \leq C_{4} {\epsilon}$$

It remains to prove that $l_{2} \leq l_{1}$.Let us assume $l_{1} < l_{2}$.Then there are two numbers $\epsilon$ and $\epsilon_{1}$,both positive such that there exists
$f_{\epsilon_{1}}$ ,so that this function is in $A^p_\nu$ and $\epsilon > \epsilon_{1}$ and also the following  conditions holds

$$\|f-f_{\epsilon_{1}}\|_{A^\infty_t}\leq \epsilon_{1}$$ and $G(f)=\infty$ ,where $G$ is a double integral  in formulation of theorem in $l_{2}$  .
(see (\ref{ETd1b}))

Next from $$\|f-f_{\epsilon_{1}}\|_{A^\infty_t} \leq \epsilon_{1}$$ we have the following two estimates,the second one is a direct corrolary of first one.First we have

for $z=x+iy$

$$(\epsilon-\epsilon_{1})\tau_{V_{\epsilon,t}}(z)\Delta^{-t}(y)\leq C |f_{\epsilon_{1}}(z)|$$,where $\tau_{V_{\epsilon,t}}(z)$ is a characteristic function of $V=V_{\epsilon,t}(f)$ set we defined above. 
 
And from last estimate we have directly multiplying both sides by Bergman kernel $B_\beta(z,w)$ and integrating by tube $T_\Omega$ both sides with measure $dV_{\beta}$

$$G(f)\leq C \int_{T_\Omega}( L(f_{\epsilon_{1}}))^{p}\Delta^{\nu-\frac{n}{r}}(y) dy dx$$, where $$L=L(f_{\epsilon_{1}},z)$$ and

$$L(f_{\epsilon_{1}},z)=\int_{T_\Omega}|f_{\epsilon_{1}}(w)||B_\beta(z,w)|dV_{\beta}(w)$$.Denote this expression by $I$.Put $\beta+\frac{n}{r}=k_{1}+k_{2}$,

where $k_{1}=\beta-\frac{n}{r}-\mu$,$k_{2}= \mu + 2\frac{n}{r}(\frac{1}{p}+\frac{1}{p_{1}})$

By classical Holder inequality  with $p$ and $p_{1}$,$p^{-1}+p_{1}^{-1}=1$ we  obviously have

$$I^{p} \leq C I_{1} I_{2}$$, where

$$I_{1}(f)= \int_{T_\Omega}|f_{1}(z)|^{p}|\Delta^{s}((z-\overline w)/i)|\Delta^{(\beta-\frac{n}{r})p}(y)dõdy$$

$$I_{2}^{\frac{p_{1}}{p}}= \int_{T_\Omega}|\Delta^{v}((z-\overline w)/i)|dxdy$$  

where $f_{1}=f_{\epsilon_{1}}$ and $$s=\mu p-2\frac{n}{r}-\beta p+p \frac{n}{r}$$,$$v=-2\frac{n}{r}-\mu p_{1}$$.

Choosing finnaly $\mu$ ,so that the estimate (A) can be used twice above and finnaly making some additional calculations we will get what we need.

Note here we  have to use the fact that $$\nu>p(\frac{n}{r}-1)$$ which was given in formulation of our theorem.

Hence we have now,

$$\int_{T_\Omega}(\int_{T_\Omega} |f_{\epsilon_{1}}(z)|B_{\beta}(z,w)|dV_{\beta}(z))^{p}\Delta^{\nu-\frac{n}{r}}(v)dV(w) \leq C \|f_{\epsilon_{1}}\|^{p}_{A^{p}_{\nu}}$$

$$G(f)\leq C \|f_{\epsilon_{1}}\|_{A^p_\nu}$$,but we also have $$f_{\epsilon_{1}} \in A^p_\nu $$.

This will give as a contradiction with our assumption above that $$G(f)=\infty$$. 
So we proved the estimate which we wanted to prove.The proof of our first theorem is now complete.

Finnaly we add some vital remarks .Similar  results are also true for certain analytic Besov spaces namely $B^{p,q}_\nu$ classes in tube domains over symmetric cones.To get such  result for these classes
we have to repeat arguments in proofs of theorems above  and use at final step in proofs the embedding theorems which connect them directly with Bergman spaces. (see \cite{DD},\cite{BBGR})
For almost all facts we mention below we refer the reader  ,for example, to \cite{BBGR}.We refer the reader for definition of analytic Besov $B^{p,q}_\nu$ spaces to \cite{BBGR}.
In the following result we use the notation $q_{\nu,p}=(\min({p,p_{1}}))q_\nu$,$q_\nu=\frac{\nu+\frac{n}{r}-1}{\frac{n}{r}-1}$,$\frac{1}{p}+\frac{1}{p_{1}}=1$
The problem on distances in analytic Besov spaces can still be posed since  
$$ B^{p,q}_\nu \subset A^\infty_{\frac{n}{rp}+\frac{\nu}{q}} $$ for $\nu>0$, $1\leq p<\infty$,$1<q < Q_{\nu,p}$,
for certain $Q_{\nu,p}$, $Q_{\nu,p}=\frac{\nu+\frac{n}{r}-1}{L}$ ,$L=\max{({0},{\frac{n}{rp_{1}}-1})}$,for $\frac{n}{r}\leq p_{1}$ we put $Q_{\nu,p}=\infty$,\cite{SE}
Note the following estimate is true $q_{\nu,p}\leq Q_{\nu,p}$.
Various projection theorems in analytic Besov spaces in tube domains over symmetric cones  were well studied recently in \cite{BBGR}.
Note $B^{2}_{0}=H^{2}$,where $H^{p}$ is analytic HArdy space in tube domain over symmetric cone.

We remark (see \cite{BBGR}) that the $A^{p,q}_\nu$ class for values 
of parametrs which we consider  is a dence subspace of $B^{p,q}_\nu$ and hence  the Bergman representaion formula is valid for all functions from this Besov class(\cite{BBGR}.)
Note also the Bergman representation formula with Bergman kernel with index $\alpha$ is valid for all functions from $A^{p,q}_\nu$ for $\alpha>\frac{n}{r}-1$ and for
$1\leq p<\infty$,and $1\leq q<Q_{\nu,p}$ .(see for example \cite{SE} and referencs there)
We formulate our first theorem via analytic Besov spaces below.
The restriction  on $p$ from formulation of previous theorems will be replaced by new restrictions on $p$ via $Q_{\nu,p}$ and $q_{\nu,p}$ which 
allows to use embeddings we need connecting analytic Besov and Bergman classes  (see \cite{SE},\cite{BBGR}) 
and also allows to pose a problem on distances in Besov class.

The following theorem follows directly from previous theorem  and the fact that for values of parametrs in our theorem below analytic Besov spaces and Bergman spaces simply coincide.(see \cite{BBGR})

\begin{thm}\label{Td1}
Let $1 < p < q_{\nu,p}$,  $\nu > p(\frac{n}{r}-1)$,$\beta > t+\frac{n}{r}-1$ ,
$t = \frac{1}{p} (\nu + \frac{n}{r})$.Set, for $f \in A^\infty_{\frac{n}{rp} + \frac{\nu}{p}}$:

\begin{equation}\label{ETd1a}
l_1(f) = {\rm dist}_{A^\infty_{\frac{n}{rp} + \frac{\nu}{p}}} (f, B^p_\nu),
\end{equation}
\begin{equation}\label{ETd1b}
l_2(f) = \inf \left\{ \epsilon > 0 : \int_{T_\Omega} \left( \int_{V_{\epsilon, t}(f)} \frac{\Delta^{\beta - t - \frac{n}{r}}(v)
dudv}{\Delta^{\beta + \frac{n}{r}}((z-\overline w)/i )} \right)^p \Delta^{\nu-\frac{n}{r}}(y) dxdy < \infty \right\}.
\end{equation}

Then there is a positive number $\beta_{0}$   ,so that for all $ \beta > \beta_{0}$   we have $l_1(f) \asymp l_2(f)$.
\end{thm}

The following theorem follows directly from the fact that for values of parametrs in our theorem below analytic Besov space and analytic Bergman classes
are connected via the following embedding  $A^{p}_\nu \subset B^{p}_\nu$ (see \cite{BBGR})

\begin{thm}\label{Td1}
Let $1 < p < Q_{\nu,p}$,  $\nu > p(\frac{n}{r}-1)$,$\beta > t+\frac{n}{r}-1$ ,
$t = \frac{1}{p} (\nu + \frac{n}{r})$.Set, for $f \in A^\infty_{\frac{n}{rp} + \frac{\nu}{p}}$:

\begin{equation}\label{ETd1a}
l_1(f) = {\rm dist}_{A^\infty_{\frac{n}{rp} + \frac{\nu}{p}}} (f, B^p_\nu),
\end{equation}
\begin{equation}\label{ETd1b}
l_2(f) = \inf \left\{ \epsilon > 0 : \int_{T_\Omega} \left( \int_{V_{\epsilon, t}(f)} \frac{\Delta^{\beta - t - \frac{n}{r}}(v)
dudv}{\Delta^{\beta + \frac{n}{r}}((z-\overline w)/i )} \right)^p \Delta^{\nu-\frac{n}{r}}(y) dxdy < \infty \right\}.
\end{equation}

Then there is a positive number $\beta_{0}$   ,so that for all $ \beta > \beta_{0}$   we have $l_1(f) < C l_2(f)$.

\end{thm}

Using obvious embeddings between $A^{p,q}_\nu$ and $A^{p}_\nu$ ,( $B^{p,q}_\nu$ and $B^{p}_\nu$) other (not sharp) similar assertions on distances for mixed norm $A^{p,q}_\nu$ and $B^{p,q}_\nu$ 
spaces as corrolaries from theorems we formulated can be also easily provided.For example,it is known $B^{p}_\nu \subset B^{q}_\nu$ for $p\leq q$.
or $A^{p,q}_\nu \subset A^{p,t}_\beta$,for all $1\leq q \leq t<\infty$ and $\frac{\nu}{q}=\frac{\beta}t$ with appropriate estimates \cite{BBR} 
The crucial fact here is that the distance problem can be posed for both spaces we discussed, since ,for example, as we mentioned above
$$ B^{p,q}_\nu \subset A^\infty_{\frac{n}{rp}+\frac{\nu}{q}}$$ for $\nu>0$, $1\leq p<\infty$,$1<q < Q_{\nu,p}$ and if we replace the BEsov class by $$A^{p,q}_\nu$$ space
it is also true ,but with other restrictions on parametrs and moreover the Bergman representation is also valid.

In one more assertion below we provide a variant of first theorem  for remaining values of positive $p$ ,assuming the embedding between two Bergman spaces $A^\infty_t$
and $A^p_\nu$ and the Bergman representation formula are still valid and hence a similar distance problem can be  posed again here.

Let $0 < p < 1$, $\nu > p(\frac{n}{r}-1)$, $\beta > t+\frac{n}{r}-1$,
$t = \frac{1}{p} (\nu + \frac{n}{r})$. Set, for $f \in A^\infty_{\frac{n}{rp} + \frac{\nu}{p}}$:
$$l_1(f) = {\rm dist}_{A^\infty_{\frac{n}{rp} + \frac{\nu}{p}}} (f, A^p_\nu)$$,

$$l_2(f) = \inf \left\{ \epsilon > 0 : \int_{T_\Omega} \left( \int_{V_{\epsilon, t}(F)} \frac{\Delta^{\beta - t - \frac{n}{r}}(v)
dudv}{\Delta^{\beta + \frac{n}{r}}((z-\overline w) /i )} \right)^p \Delta^{\nu-\frac{n}{r}}(y) dxdy < \infty \right\}.$$

Assuming that the embedding $A^p_\nu \subset  A^\infty_t$ is true and the Bergman  representation via $B_\beta(z,w)$ kernel holds 
there is a positive number $\beta_{0}$ , so that for all $\beta > \beta_{0}$ 
we have $l_1(f) \leq C l_2(f)$,where $C$ is a positive constant depending on parametrs.
The proof of this theorem follows directly from the proof of first theorem ,namely this proof ,as it can be easily seen ,coincides with the proof of first part of first theorem.   
We finnaly formulate one more theorem.
Let $H^{p}(T_\Omega)$ be standard Hardy space in tubular domain over symmetric cone.Let further $H^{p}_\alpha(T_\Omega)$ be weighted Hardy class with $\Delta^{\alpha}(y)$ weight.\cite{BBR}

These are spaces with finite norm $$\sup_{y \in \Omega}\int_{R^{n}}|f(\tau+iy)|^{p}d\tau \Delta^{p\alpha}(y)$$,$1\leq p <\infty$,$\alpha \in R$

Let $1\leq p< \infty$,$1\leq q < \infty$.Let also $q\leq s$.Then (see \cite{BBR}) $A^{p,q}_\nu \subset H^{s}_\beta$ ,where $\nu>\frac{n}{r}-1$ and 
where $\beta=\frac{\nu}{q}+\frac{n}{rp}-\frac{n}{rs}$ .For $p=q=s$ this embedding with appropriate estimate is taking obviously  a very simple
form (see \cite{BBR}(proposition 3.5)) and the distance problem here can be easily posed again obviously in general case and in mentioned simple case.
Note for analytic and harmonic
function spaces it was posed and solved in \cite{AS1} and \cite{SM1} and \cite{SM2}.Repeating arguments from these papers an analogue of that result using approaches we provided above can also
be obtained.In particular the following theorem is true.

Let
\begin{equation}\label{Efsep}
L_{\epsilon, s}(f) = \left\{ y \in \Omega : \int_{R^{n}} |f(\tau+iy)| d\tau \Delta^s(y) \geq \epsilon \right\}
\end{equation}
for fixed positive $\epsilon$ and analytic in tube $f$ function.
By $\tau_M$ as before we denote a characteristic function of $M$ set. 

\begin{thm}\label{Td1}

Let $\nu > \frac{n}{r}-1$.
Set, for $f \in H^{1}_{\nu}$:
\begin{equation}\label{ETd1a}
l_1(f) = {\rm dist}_{H^{1}_{\nu}}(f, A^1_\nu),
\end{equation}
\begin{equation}\label{ETd1b}
l_2(f) = \inf \left\{ \epsilon > 0 : \int_{\Omega} \left( \tau_{L_{\epsilon, t}(f)}(y)  \right) \Delta^{-\frac{n}{r}}(y) dy < \infty \right\}.
\end{equation}

Then   $l_1(f) \asymp l_2(f)$.

\end{thm}

The proof of one implication $l_{2}(f) \leq l_{1}(f)$ of this sharp result can be obtained easily by small modification of arguments from proof of first theorem above 
or by repetition of arguments from our proof of parallel mentioned result in one dimension from our previous just mentioned papers.The other implication is obtained by us with the help 
of approaches we provided in proof of our first theorem again following arguments of onedimensional result \cite{SA} and \cite{AS1}  and estimates we used above for Bergman kernel
and $\Delta(y)$ function.We leave these detailes to readers noting that the complete proof will be presented elsewhere.

\end{document}